\documentclass[12pt]{amsart}
\usepackage{times}
\usepackage{amssymb}
\usepackage{amsfonts}
\usepackage{amsthm}
\usepackage{amsmath}
\usepackage{comment}
\advance \topmargin by -\headheight
\advance \topmargin by -\headsep
\evensidemargin \oddsidemargin
\DeclareMathOperator{\DT}{DT}
\DeclareMathOperator{\tr}{tr}
\newcommand\Bc{{\mathcal{B}}}
\newcommand\Cpx{{\mathbb C}}
\newcommand\eps{\epsilon}
\newcommand\Fb{{\mathbb F}}
\newcommand\ImagPart{{\operatorname{Im}}}
\newcommand\Mcal{{\mathcal{M}}}
\newcommand\Mcalt{{\widetilde\Mcal}}
\newcommand\Nats{{\mathbb N}}
\newcommand\Nc{{\mathcal{N}}} 
\newcommand\Oc{{\mathcal{O}}} 
\newcommand\RealPart{{\operatorname{Re}}}
\newcommand\restrict{{\upharpoonright}}
\newcommand\taut{{\tilde\tau}}
\newcommand\vol{{\operatorname{vol}}}
\newcommand\Yt{{\widetilde Y}}
\newcommand\Zt{{\widetilde Z}}
\newcommand\zt{{\widetilde z}}

\newcount\theTime
\newcount\theHour
\newcount\theMinute
\newcount\theMinuteTens
\newcount\theScratch
\theTime=\number\time
\theHour=\theTime
\divide\theHour by 60
\theScratch=\theHour
\multiply\theScratch by 60
\theMinute=\theTime
\advance\theMinute by -\theScratch
\theMinuteTens=\theMinute
\divide\theMinuteTens by 10
\theScratch=\theMinuteTens
\multiply\theScratch by 10
\advance\theMinute by -\theScratch

\def\today{{\number\day\space
 \ifcase\month\or
  January\or February\or March\or April\or May\or June\or
  July\or August\or September\or October\or November\or December\fi
 \space\number\year}}

\newtheorem{theorem}{Theorem}[section]
\newtheorem{lemma}[theorem]{Lemma}

\theoremstyle{remark}

\begin{document}

\title[Free Dimension of $\DT$--operators]
{The Microstates Free Entropy Dimension of any DT--operator is 2}

\date{14 December 2004}
\author[K. Dykema, K. Jung, D. Shlyakhtenko]
{Ken Dykema*, Kenley Jung$\dagger$, and Dimitri Shlyakhtenko$\ddagger$}

\address{K. Dykema, Mathematisches Institut, Westf\"alische Wilhelms--Universit\"at M\"unster,
Einsteinstr.\ 62, 48149 M\"unster, Germany;
{\rm permanent address:} Department of Mathematics, Texas A\&M University,
College Station, TX 77843-3368, USA}
\email{kdykema@math.tamu.edu}

\address{K.\ Jung, D.\ Shlyakhtenko, Department of Mathematics, University of California,
Los Angeles, CA 90095-3840, USA}
\email{kjung@math.ucla.edu, shlyakht@math.ucla.edu}

\thanks{*Research supported by the Alexander von Humboldt Foundation
and NSF grant DMS--0300336. \\
$\dagger$ Research supported by an NSF Postdoctoral Fellowship. \\
$\ddagger$ Research supported by a Sloan Foundation fellowship and 
NSF grant DMS-0355226.}

\begin{abstract} Suppose that $\mu$ is an arbitrary Borel measure on
$\mathbb C$ with compact support and $c >0$.  If $Z$ is a $\DT(\mu,
c)$-operator as defined by Dykema and Haagerup in 
\cite{dykema-haagerup:DT}, then the microstates free entropy dimension of $Z$ is $2$.
\end{abstract}
\maketitle

\section{Introduction.}

$\DT$--operators were introduced by Dykema and Haagerup in their work on 
invariant subspaces of certain operators in a II$_1$ factor
\cite{dykema-haagerup:invariant,dykema-haagerup:DT}.  A $\DT$--operator
$Z$
is specified by two parameters, $\mu$ and $c$, where $c>0$ and 
$\mu$ is a Borel probability measure on $\mathbb C$ with compact
support.  Roughly, the operator $Z$ is determined by stating that its
$*$--distribution is the same as the limit $*$--distribution 
as $N\to\infty$ of random matrices
$$Z_N = D_N + c T_N,$$
where $D_N$ are diagonal $N\times N$ matrices whose spectral measures
converge to $\mu$ in distribution, while $T_N$ is a 
strictly upper triangular random $N\times N$ matrix with i.i.d.\ Gaussian
entries.    
Equivalently, (see \cite{Sn}, \cite{shlyakht:bandmatrix}, \cite{dykema-haagerup:DT} and the appendix of~\cite{DH3}),
$Z$ can be viewed as a  sum $ Z = d + c T$, where $d$ is a normal 
operator with spectral measure $\mu$
contained in a diffuse von Neumann algebra $A$, and $T$ is an 
$A$-valued circular operator with a certain covariance.
Finally, a result of \'Sniady~\cite{sniady} shows that 
a $\DT(\mu,c)$--operator is one whose free entropy is maximized among
all those operators having Brown measure equal to $\mu$ and with a fixed off--diagonality.

If we write $Z=d+cT$ as above, it is clear that $W^*(Z)\subset W^*(d,T)\subseteq W^*(A\cup\{T\})$,
while a simple computation shows $W^*(A\cup\{T\})=L(\mathbb{F}_2)$.
By Lemma~6.2 of ~\cite{dykema-haagerup:DT}, for any $\mu$
we may choose $d$ having trace of spectral measure equal to $\mu$ and so
that  $d,T\in W^*(Z)$;
by~\cite{DH3}, $A\subseteq W^*(T)$,
so we always have $W^*(Z)\cong L(\mathbb{F}_2)$.
Thus $Z$ can be viewed as an interesting generator for this free group
factor.

In order to test the hypothesis that Voiculescu's free entropy dimension
$\delta_0$ 
\cite{dvv:entropy2,dvv:entropy3,dvv:survey} is the same for any sets
of generators of a von Neumann algebra, it is important to decide whether
the free entropy dimension of $Z$ is $2$ ($L(\mathbb{F}_2)$ 
clearly has another set of generators of free entropy dimension $2$).

For another version of free entropy dimension, also defined by Voiculescu, called the non-microstates
free entropy dimension \cite{dvv:entropy5}, L. Aagaard has recently
shown \cite{aagard:dt} that the dimension of $Z$ is indeed $2$.
It is known by~\cite{BCG}
that the non-microstates free entropy dimension dominates $\delta_0$
but at present it is open whether the reverse inequality holds.
Thus, Aagaard's result does not solve the question for the original 
microstates definition.

In this paper, we show that, indeed, $\delta_0(Z) = 2$. Our proof uses an
equivalent packing number formulation of the microstates free entropy dimension,
due to Jung~\cite{jung:freeEntropyLemma}. In this approach, to get the nontrivial
lower bound on $\delta_0(Z)$, one must have lower bounds on the
$\epsilon$--packing numbers of spaces of matricial microstates for $Z$, which are
in turn obtained by lower bounds on the volume of $\epsilon$--neighborhoods of
these microstate spaces. The $k$th microstate space is the set
$\Gamma(Z;m,k,\gamma)$, for $m,k\in\Nats$ and $\gamma>0$, of all $k\times k$
complex matrices whose $*$--moments up to order $m$ are $\gamma$--close to the
values of the corresponding $*$--moments of $Z$, and the volumes are for Lebesgue
measure $\lambda_k$ on $M_k(\Cpx)$ viewed as a Euclidean space of real dimension
$2k^2$ with coordinates corresponding to the real and imaginary parts of the
entries of a matrix.

In order to outline how we get these lower bounds on volumes, let us for
convenience take $Z$ equal to the $\DT(\delta_0,1)$--operator $T$. A key result
that we use is a recent one of Aagaard and Haagerup~\cite{aagard-haagerup:dt},
showing that a certain $\eps$--perturbation of $T$ has Brown measure uniformly
distributed on the disk of radius $r_\eps:=1/\sqrt{\log(1+\eps^{-2})}$ centered
at the origin; note how slowly this disk shrinks as $\eps$ approaches zero.
Applying a result of \'Sniady~\cite{Sn1} to this situation, we find matrices
$A_k\in M_k(\Cpx)$ that lie in $\eps$--neighborhoods of microstate spaces for
$T$, whose eigenvalues are close to uniformly distributed (as $k$ gets large) in
the disk of radius $r_\eps$. Thus, in order to get a lower bound on the volume of
a $2\eps$--neighborhood of a microstate space for $T$, it will suffice to get a
lower bound on the volume of a unitary orbit of an $\eps$--neighborhood of $A_k$.

Every element of $M_k(\Cpx)$ has an upper triangular matrix in its unitary orbit.
Thus, letting $T_k(\Cpx)$ denote the set of upper triangular matrices in
$M_k(\Cpx)$, there is a measure $\nu_k$ on $T_k(\Cpx)$ such that
$\lambda_k(\Oc)=\nu_k(\Oc\cap T_k)$ for every $\Oc\subseteq M_k(\Cpx)$
invariant under unitary conjugation.
Freeman Dyson identified such a measure $\nu_k$ (see Appendix~35
of~\cite{Mehta}), and showed that if we view $T_k(\Cpx)$ as a Euclidean space of
real dimension $k(k-1)$ with coordinates corresponding to the real and imaginary
parts of the matrix entries lying on and above the diagonal, then $\nu_k$ is
absolutely continuous with respect to Lebesgue measure on $T_k(\Cpx)$ and has
density given at $B=(b_{ij})_{1\le i,j\le k}\in T_k(\Cpx)$ by
\begin{equation}\label{eq:DysonDensity} C_k\prod_{1\le p<q\le
k}|b_{pp}-b_{qq}|^2, \end{equation} where the constant is
\begin{equation}\label{eq:Ck} C_k =\frac{\pi^{k(k-1)/2}}{\prod_{j=1}^{k}j!}.
\end{equation}

We will use this measure of Dyson to find lower bound on the volume of unitary
orbits of an $\eps$--neighborhood of $A_k$, and we may take $A_k$ to be upper
triangular. However, so far we only have information about the eigenvalues of
$A_k$, namely the diagonal part of it. Loosely speaking, in order to get a handle
on the part strictly above the diagonal, we use a result of Dykema and
Haagerup~\cite{dykema-haagerup:DT} to realize $T$ as an upper triangular matrix
\[ T= \frac1{\sqrt N} \begin{bmatrix} T_{11} & T_{12} & \cdots & T_{1N} \\ 0 &
T_{22} & \ddots & \vdots \\ \vdots & \ddots & \ddots & T_{N-1,N} \\ 0 & \cdots &
0 & T_{NN} \end{bmatrix} \] of operators where each $T_{ii}$ is a copy of $T$,
each $T_{ij}$ for $i<j$ is circular and the family $(T_{ij})_{1\le i\le j\le N}$
is $*$--free. Thus, $A_k$ can be taken to be of the form \[ \begin{bmatrix}
B_{11} & B_{12} & \cdots & B_{1N} \\ 0 & B_{22} & \ddots & \vdots \\ \vdots &
\ddots & \ddots & B_{N-1,N} \\ 0 & \cdots & 0 & B_{NN} \end{bmatrix} \] where
each $B_{ii}$ is upper triangular, where we have good knowledge of the eigenvalue
distributions of each $B_{ii}$ and where the $B_{ij}$ for $i<j$ approximate
$*$--free circular elements.
Using the strengthened asymptotic freeness results of Voiculescu~\cite{dvv:strenghthened},
we find enough approximants for these $B_{ij}$.
Although we still have no real knowledge about the
entries of the $B_{ii}$ lying above the diagonal, these parts are of negligibly
small dimension as $N$ gets large, and we are able to get good enough lower
bounds. The techniques we use for estimating integrals of the quantity~\eqref{eq:DysonDensity}
over certain regions are taken from~\cite{jung:fractalDimension}.

\section{Microstates for $Z$ with well--spaced spectral densities}

The following lemma is an application of the result of Aagaard and Haagerup~\cite{aagard-haagerup:dt}
mentioned in the introduction in order to make perturbations of general DT--operators having Brown
measure that is relatively well spread out.
For an element $a$ of a noncommutative probability space $(\Mcal,\tau)$, we write $\|a\|_2$ for $\tau(a^*a)^{1/2}$.

\begin{lemma}\label{lem:Brown}
Let $\mu$ be a compactly supported Borel probability measure on $\Cpx$ and let $c>0$.
Let $Z$ be a $\DT(\mu,c)$--operator in a W$^*$--noncommutative probability space
$(\Mcal,\tau)$.
Let us write
\[
\mu=\nu+\sum_{i=1}^sa_i\delta_{z_i}
\]
for some $s\in\{0\}\cup\Nats\cup\{\infty\}$, $z_i\in\Cpx$ and $a_i>0$, where $\nu$ is a diffuse
measure and where $z_i\ne z_j$ if $i\ne j$.
Consider the W$^*$--noncommutative probability space
\[
(\Mcalt,\taut)=(\Mcal,\tau)*(L(\Fb_2),\tau_{\Fb_2}).
\]
Then for every $\eps>0$, there is $\Zt_\eps\in\Mcalt$ such that $\|\Zt_\eps-Z\|_2\le\eps c$ and
where the Brown measure of $\Zt_\eps$ is equal to
\[
\sigma_\eps:=\nu+\sum_{i=1}^sa_i\rho_{i,\eps},
\]
where $\rho_{i,\eps}$ is the probability measure that is uniform distribution on the disk centered
at $z_i$ and having radius
\[
r_i:=c\sqrt{\frac{a_i}{\log(1+a_i\eps^{-2})}}.
\]
Finally, if $\delta>0$ and if
\[
X_\delta=\{(w_1,w_2)\in\Cpx^2\mid|w_1-w_2|<\delta\},
\]
then
\begin{equation}\label{eq:sigmaX}
(\sigma_\eps\times\sigma_\eps)(X_\delta)\le(\nu\times\nu)(X_\delta)
+2\sum_{i=1}^s\min(a_i,\delta^2c^{-2}\log\big(1+a_i\eps^{-2})).
\end{equation}
\end{lemma}
\begin{proof}
By results from~\cite{dykema-haagerup:DT}, taking projections onto local spectral subspaces
of $Z$, we find projections $p_j\in\Mcal$ (for $0\le j<s+1$) such that
\begin{itemize}
\item $\sum_{j=0}^sp_j=1$,
\item $p_0+p_1+\cdots+p_k$ is $Z$--invariant for all integers $k$ such that $0\le k<s+1$,
\item $\tau(p_k)=\begin{cases}|\nu|&\text{if }k=0 \\ a_k&\text{if }1\le k<s+1,\end{cases}$
\item In $(p_k\Mcal p_k,\tau(p_k)^{-1}\tau\restrict_{p_k\Mcal p_k})$, $p_kZp_k$ is $\DT(|\nu|^{-1}\nu,c\sqrt{|\nu|})$
  if $k=0$ and is $\DT(\delta_{z_k},c\sqrt{a_k})$ if $1\le k<s+1$.
\end{itemize}
Let $Y\in\Mcalt$ be centered circular such that $Y$ and $Z$ are $*$--free and $\taut(Y^*Y)=1$.
Let
\begin{equation}\label{eq:ZY}
\Zt_\eps=Z+\eps\sum_{i=1}^sa_i^{-1/2}cp_iYp_i.
\end{equation}
Then $\|\Zt_\eps-Z\|_2^2=\eps^2c^2\sum_{i=1}^sa_i\le\eps^2c^2$.
On the other hand, $\Zt_\eps$ is upper triangular with respect to the projections $p_0,p_1,\ldots$;
the Brown measure of $\Zt_\eps$ is, therefore, equal to the Brown measure of its diagonal part
\begin{equation}\label{eq:Ztdiag}
p_0Zp_0+\sum_{i=1}^s\big(p_iZp_i+\eps\,a_i^{-1/2}cp_iYp_i\big).
\end{equation}
But in $(p_i\Mcalt p_i,a_i^{-1}\taut\restrict_{p_i\Mcalt p_i})$, the operator $\eps\,a_i^{-1/2}cp_iYp_i$ is a 
centered circular operator of second moment $\eps^2c^2$
that is $*$--free from the $\DT(\delta_{z_i},c\sqrt{a_i})$ operator $p_iZp_i$.
Therefore, the random variable
\begin{equation}\label{eq:piZ}
p_iZp_i+\eps\,a_i^{-1/2}cp_iYp_i
\end{equation}
has the same $*$--distribution as
$z_iI+c\sqrt{a_i}(T+\eps\,a_i^{-1/2}Y)$, where $T$ is a $\DT(\delta_0,1)$--operator that is $*$--free from $Y$.
By~\cite{aagard-haagerup:dt}, the Brown measure of the random variable~\eqref{eq:piZ} is equal to $\rho_{i,\eps}$.
This yields $\sigma_\eps$ for the Brown measure of the operator~\eqref{eq:Ztdiag}, hence of $\Zt_\eps$ itself.

Finally, we have
\begin{equation}\label{eq:sigsig}
(\sigma_\eps\times\sigma_\eps)(X_\delta)\le(\nu\times\nu)(X_\delta)
+2\sum_{i=1}^sa_i(\sigma_\eps\times\rho_{i,\eps})(X_\delta)
\end{equation}
and
\begin{equation}\label{eq:sigrho}
(\sigma_\eps\times\rho_{i,\eps})(X_\delta)=\int_\Cpx\rho_{i,\eps}(w+\delta\,{\mathbb D})d\sigma_\eps(w)
\le\min(1,\delta^2r_i^{-2}),
\end{equation}
where $\mathbb D$ is the unit disk in $\Cpx$.
Taken together,~\eqref{eq:sigsig} and~\eqref{eq:sigrho}
yield the inequality~\eqref{eq:sigmaX}.
\end{proof}

The next lemma uses a result of \'Sniady~\cite{Sn1} to find matrix
approximants of the operators appearing in Lemma~\ref{lem:Brown}.

In the following lemma and throughout this paper, for a matrix $A\in M_k(\Cpx)$ we let
$|A|_2=\tr_k(A^*A)^{1/2}$, where $\tr_k$ is the normalized trace on $M_k(\Cpx)$.
Moreover, by the eigenvalue distribution of $A\in M_k(\Cpx)$
we mean its Brown measure, which is just the probability
measure that is uniformly distributed on its list of eigenvalues $\lambda_1,\ldots,\lambda_k$,
where these are 
listed according to (general) multiplicity,
i.e.\  a value $z$ is listed $\dim\bigcup_{n=1}^\infty\ker((A-zI)^n)$ times.

\begin{lemma}\label{lem:ms}
Let $\mu$ be a compactly supported Borel probability measure on $\Cpx$ and let $c>0$.
Then there 
exists a sequence $\langle y_k \rangle_{k=1}^{\infty}$ such that for any $\eps>0$,
there exists a sequence $\langle z_{k, \epsilon} \rangle_{k=1}^\infty$ such that 
\begin{itemize}
\item $y_k, z_{k,\epsilon}\in M_k(\Cpx)$,
\item $\|y_k\|$ and $\|z_{k, \epsilon}\|$ remain bounded as $k\to\infty$,
\item $\limsup_{k\to\infty}|y_k-z_{k,\epsilon}|_2\le\eps c$,
\item $y_k$ converges in $*$--moments as $k\to\infty$ to a 
$\DT(\mu,c)$--operator,
\item the eigenvalue distribution of $z_{k,\epsilon}$ converges weakly as 
$k\to\infty$ to 
the measure
$\sigma_\eps$ described in Lemma~\ref{lem:Brown}.
\end{itemize}
\end{lemma}
\begin{proof}
Let $Z$ be a $\DT(\mu,c)$--operator, let $\Yt$ be the operator
$\sum_{i=1}^sa_i^{-1/2}cp_iYp_i$
appearing in~\eqref{eq:ZY} in the proof of the preceding lemma, so that $\Zt_\eps=Z+\eps\Yt$.
Since $Z$ can be constructed in $L(\Fb_2)$ and since free group factors
can be embedded in the ultrapower $R^\omega$ of the hyperfinite II$_1$ factor, 
there are bounded
sequences $\langle y_k \rangle_{k=1}^\infty$ and $\langle d_k 
\rangle_{k=1}^\infty$ such that $y_k,d_k\in M_k(\Cpx)$ and such that the
pair $y_k,d_k$ converges in $*$--moments to the pair $Z,\Yt$.
Letting $\zt_k=y_k+\eps d_k$, we have that $\zt_k$ converges in $*$--moments 
to $\Zt_\eps$ as $k\to\infty$.
By Theorem~7 of~\cite{Sn1}, there is a sequence $\langle z_{k, \epsilon} 
\rangle_{k=1}^\infty$ with $z_{k, \epsilon} \in M_k(\Cpx)$ such that
$\|z_{k, \epsilon} -\zt_{k, \epsilon}\|$ tends to zero and the eigenvalue 
distribution of $z_{k, \epsilon}$ converges weakly as $k\to\infty$
to the Brown measure of $\Zt_\eps$, namely, to $\sigma_\eps$.
\end{proof}

Suppose that $\lambda = \langle \lambda_j \rangle_{j=1}^k$ is a finite sequence
of complex numbers.
For each $j$, write $\lambda_j = a_j + i b_j$, $a_j, b_j \in
\mathbb R.$ Define $ Q_{\epsilon} = \prod_{j=1}^k [a_j - \epsilon, a_j +
\epsilon]$ and $R_{\epsilon} = \prod_{j=1}^k [b_j -\epsilon, b_j + \epsilon]$.
Set
\[
E_\eps(\lambda) = \int_{R_{\epsilon}} \bigg( \int_{Q_{\epsilon}}
\prod_{\substack{1 \leq i,j \leq k \\ i\ne j}} \left( |s_i - s_j| + | t_i - t_j|^2 \right )^{1/2}
ds \bigg )  dt,
\]
where $ds = ds_1 \cdots ds_k$ and $dt = dt_1 \cdots dt_k$.

The following lemma proves lower bounds for certain asymptotics of the quantities $E_\eps(\lambda)$.
We will apply this lemma to the case when $\lambda$ is the eigenvalue sequence
of matrices like the $z_{k,\eps}$ found in Lemma~\ref{lem:ms}.

\begin{lemma}\label{lem:Eeps}
Let $\mu$ and $c$ be as in Lemma~\ref{lem:Brown}.
For each $\epsilon >0$ and $k\in\Nats$, let
$\lambda^{(k,\eps)}=\langle\lambda^{(k,\eps)}_1,\ldots,\lambda^{(k,\eps)}_{n(k)}\rangle$
be a finite sequence of complex numbers and assume that for every $\eps>0$,
\[
\sup_{k\in\Nats,\,1\le j\le n(k)}|\lambda^{(k,\eps)}_j|<\infty
\]
and the probability measures
\begin{equation}\label{eq:kepsmeas}
\frac1{n(k)}\sum_{j=1}^{n(k)}\delta_{\lambda^{(k,\eps)}_j}
\end{equation}
converge weakly to the measure $\sigma_\eps$ of Lemma~\ref{lem:Brown} as $k\to\infty$.
Let
\[
f(\epsilon) = \liminf_{k \rightarrow \infty}n(k)^{-2}\log(E_{\epsilon}(\lambda^{(k,\eps)})).
\]
Then
\begin{equation}\label{eq:lili} \liminf_{\eps\to0}
\left( \frac{f(\epsilon)}{|\log\eps|} \right) \geq 0. \end{equation}
\end{lemma}
\begin{proof}
Note that we must have $n(k)\to\infty$ as $k\to\infty$.
Given $\epsilon >0$ small,
take $1\ge\delta>3\eps$.
Define
\[
W_{k,\eps}=\{(i,j)\in\{1,\ldots,n(k)\}^2\mid i\ne j,\, |\lambda_i^{(k,\eps)}-\lambda_j^{(k,\eps)}|<\delta\}.
\]
Writing for
each $1 \leq j \leq k$, $\lambda^{(k,\eps)}_j = a_j + i b_j$ where $a_j, b_j \in \mathbb R$ define
$Q_{\epsilon, k} = \prod_{j=1}^{n(k)} [a_j - \epsilon, a_j + \epsilon]$, $R_{\epsilon,k} = \prod_{j=1}^{n(k)}
[b_j - \epsilon, b_j +\epsilon]$, and $K_{\epsilon,k} = Q_{\epsilon, k} \times R_{\epsilon,k}.$
Now
\begin{align*}
E_\eps(\lambda^{(k,\eps)})&=\int_{K_{\eps,k}}\prod_{i\ne
j}(|s_i-s_j|^2+|t_i-t_j|^2)^{1/2}dsdt \\
&\ge(\delta-\sqrt8\eps)^{n(k)^2-\#W_{k,\eps}}\int_{K_{\eps,k}}\prod_{(i,j)\in
W_{k,\eps}}(|s_i-s_j|^2+|t_i-t_j|^2)^{1/2}dsdt \\
&\ge(\delta-3\eps)^{n(k)^2-\#W_{k,\eps}}\bigg(\int_{Q_{\eps,k}}\prod_{(i,j)\in W_{k,\eps}}|s_i-s_j|ds\bigg)
\bigg(\int_{R_{\eps,k}}\prod_{(i,j)\in W_{k,\eps}}|t_i-t_j|dt\bigg),
\end{align*}
where $ds=ds_1\cdots ds_{n(k)}$ and $dt=dt_1\cdots dt_{n(k)}$.

We now wish to find a lower bounds for the two integrals in the above expression. By
Fubini's Theorem we can assume $a_1\le a_2\le\cdots\le a_{n(k)}$. Let
\[
[-\eps,\eps]^{n(k)}_<=\{(x_1,\ldots,x_{n(k)})\in[-\eps,\eps]^{n(k)}\mid x_1<x_2<\cdots<x_{n(k)}\}.
\]
Then by the change of variables
$[-\epsilon,\epsilon]^{n(k)}_{<}\ni(x_1,\ldots,x_{n(k)})\mapsto(a_1+x_1,\ldots,a_{n(k)}+x_{n(k)}) \in Q_{\epsilon,k}$
and Selberg's Integral Formula it follows that

\begin{align*} \int_{Q_{\eps,k}}\prod_{(i,j)\in W_{k,\eps}}|s_i-s_j|ds
&\ge\int_{[-\eps,\eps]^{n(k)}_<}\;\prod_{(i,j)\in W_{k,\eps}}|x_i-x_j|dx_1\cdots dx_{n(k)} \\
&\ge(2\eps)^{-(n(k)^2-n(k)-\#W_{k,\eps})} \cdot \int_{[-\eps,\eps]^{n(k)}_<}\;\prod_{i\ne
j}|x_i-x_j|dx_1\cdots dx_{n(k)} \\ &=\frac{(2\eps)^{-(n(k)^2-n(k)- \#W_{k,\eps})}}{n(k)!} \cdot
\int_{[-\eps,\eps]^{n(k)}}\prod_{i\ne j}|x_i-x_j|dx_1\cdots dx_{n(k)} \\
&=\frac{(2\eps)^{n(k)+\#W_{k,\eps}}}{n(k)!} \cdot
\prod_{j=0}^{n(k)-1}\frac{\Gamma(j+2)\Gamma(j+1)^2}{\Gamma(n(k)+j+1)}, \end{align*}

The same lower bound applies to $\int_{R_{\eps,k}}\prod_{(i,j)\in W_{k,\eps}}|t_i-t_j|dt$
so that combining these two
we get
 \begin{align*} E_\eps(\lambda^{(k,\eps)})&\ge(\delta-3\eps)^{n(k)^2-\#W_{k,\eps}}\bigg(
\frac{(2\eps)^{n(k)+\#W_{k,\eps}}}{n(k)!} 
\cdot \prod_{j=0}^{n(k)-1}\frac{\Gamma(j+2)\Gamma(j+1)^2}{\Gamma(n(k)+j+1)}\bigg)^2
\\ &\ge(\delta-3\eps)^{n(k)^2}\bigg(
\frac{(2\eps)^{n(k)+\#W_{k,\eps}}}{n(k)!} 
\cdot \prod_{j=0}^{n(k)-1}\frac{\Gamma(j+2)\Gamma(j+1)^2}{\Gamma(n(k)+j+1)}\bigg)^2.  
\end{align*} 

Using \[
\lim_{k\to\infty}n(k)^{-2}\log(\prod_{j=0}^{n(k)-1}\frac{\Gamma(j+2)\Gamma(j+1)^2}{\Gamma(n(k)+j+1)})=-2\log2,
\]
we find
\[
f(\epsilon)\ge\log(\delta-3\eps)+ 
2\log(2\eps)\,\limsup_{k\to\infty}\frac{\#W_{k,\eps}}{n(k)^2}-4\log2.
\]
Since the measures~\eqref{eq:kepsmeas} converge weakly to $\sigma_\eps$, 
by standard
approximation techniques one sees \[
\lim_{k\to\infty}\frac{\#W_{k,\eps}}{n(k)^2}=(\sigma_\eps\times\sigma_\eps)(X_\delta), \] where
$X_\delta$ is as in Lemma~\ref{lem:Brown}. As $\eps\to0$ choose $\delta=\frac1{|\log\eps|}$, so
that $\delta^2\log(1+a\eps^{-2})\to0$ for all $a>0$ and $\frac\eps\delta\to0$ and
$\frac{\log\delta}{\log\eps}\to0$. Using the upper bound~\eqref{eq:sigmaX} and the fact that
$\nu$ is diffuse, we get \[ 
\lim_{\eps\to0}(\sigma_\eps\times\sigma_\eps)(X_\delta)=0. \] Now
one easily verifies that~\eqref{eq:lili} holds. \end{proof}

\section{The Main Result}

Before beginning the main result first a few comments 
on a packing formulation for microstates free entropy dimension are in order.
If $X = 
\{x_1,\ldots,x_n\}$ is an $n$-tuple of selfadjoint elements in a tracial von 
Neumann algebra, then the free entropy dimension (as defined 
by Voiculescu~\cite{dvv:entropy3}) is given by the formula

\[ \delta_0(X) = n + \limsup_{\epsilon \rightarrow 0} \frac{\chi(x_1 + \epsilon 
s_1, \ldots, x_n + \epsilon s_n:s_1, \ldots,s_n)}{|\log \epsilon|} 
\]

\noindent where $\{s_1, \ldots, s_n\}$ is a semicircular family free from $X$.  
The packing formulation found in~\cite{jung:freeEntropyLemma} and modified
slightly in~\cite{J3} (to remove the norm restriction on microstates),
is
\[ \delta_0(X) = \limsup_{\epsilon \rightarrow 0} \frac{\mathbb 
P_{\epsilon}(X)}{|\log \epsilon|}, \]
where
\begin{equation}\label{eq:Peps}
\mathbb P_{\epsilon}(X)=\inf_{m\in\Nats,\,\gamma>0}
\limsup_{k\to\infty}k^{-2}\log P_\eps(\Gamma(X;m,k,\gamma)).
\end{equation}
Here, $\Gamma(X;m,k,\gamma)\subseteq (M_k(\Cpx)_{s.a.})^n$ is the microstate space
of Voiculescu~\cite{dvv:entropy2}, but taken without norm restriction,
as considered in~\cite{BB}, and $P_\eps$ is the packing number with respect to the metric arising
from the normalized trace.

Let $Y = \{y_1, \ldots, y_n\}$ be an arbitrary $n$-tuple 
of (possibly nonselfadjoint) elements in a tracial von Neumann algebra.
Now the definition of $\mathbb P_{\epsilon}$ makes perfect sense for the 
set $Y$ if we replace the microstate space in~\eqref{eq:Peps} 
with the non-selfadjoint $*$-microstate space
$\Gamma(Y;m,k,\gamma)\subseteq(M_k(\Cpx))^n$,
which is the set of all $n$--tuples of $k\times k$ matrices whose $*$--moments up to order $m$
approximate those of $Y$ within tolerance of $\gamma$.
Let us (temporarily) denote the
quantity so obtained by $\overline{\mathbb P_{\epsilon}}(Y)$ and define
\begin{equation}\label{eq:d0bar}
\overline{\delta_0}(Y) = \limsup_{\epsilon \rightarrow 0} \frac{
\overline{\mathbb P_{\epsilon}}(Y)}{|\log \epsilon|}.
\end{equation}
It is easy to see that if
$X$ is a set of selfadjoints, then $\overline{\mathbb P_{\epsilon}}(X) \geq 
\mathbb P_{\epsilon}(X) \geq \overline{\mathbb P_{2\epsilon}}(X)$ and that in the 
nonselfadjoint setting the quantity~\eqref{eq:d0bar} is a $*$-algebraic
invariant, so that

\begin{multline*}
\delta_0(\RealPart(y_1), \ImagPart(y_1), \ldots, \RealPart(y_n), \ImagPart(y_n))
=\limsup_{\epsilon \rightarrow 0}
  \frac{\mathbb P_{\epsilon}(\RealPart(y_1), \ImagPart(y_1), \ldots, \RealPart(y_n), \ImagPart(y_n))}
    {|\log \epsilon|}  =  \\
=\limsup_{\epsilon \rightarrow 0} \frac{\overline{\mathbb
P_{\epsilon}}(\RealPart(y_1), \ImagPart(y_1), \ldots, \RealPart(y_n), \ImagPart(y_n))}{|\log \epsilon|} = 
\limsup_{\epsilon \rightarrow 0} \frac{\overline{ \mathbb P_{\epsilon}}(Y)}{|\log\epsilon|}=\overline{\delta_0}(Y),
\end{multline*}
where $\RealPart(y_i)$ and $\ImagPart(y_i)$ are the real and imaginary parts of $y_i$.
Moreover, if $X$ is set 
of selfadjoints, then 
\[ \delta_0(X) = \limsup_{\epsilon 
\rightarrow 0} \frac{\mathbb P_{\epsilon}(X)}{|\log \epsilon|} = 
\limsup_{\epsilon \rightarrow 0} \frac{\overline{\mathbb P_{\epsilon}}(X)}{|\log 
\epsilon|}=\overline{\delta_0}(X).\]  

\noindent The following notational conventions,
which will be used in the remainder of this paper, are, therefore, justified:
for any finite set of operators $Y$
(selfadjoint or otherwise) in a tracial von Neumann algebra we will write
$\mathbb P_{\epsilon}(Y)$ for the packing quantity derived from the nonselfadjoint microstates
(that was denoted $\overline{\mathbb P_\eps}(Y)$ above)
and we will write $\delta_0(Y)$ for the
free entropy dimension of $Y$ that was denoted $\overline{\delta_0}(Y)$ above.

In the proof of the main result, we
will use $E_\eps(A)$ for $A\in M_k(\Cpx)$ to mean $E_\eps(\lambda)$,
where $\lambda = \langle \lambda_j \rangle_{j=1}^k$ are the eigenvalues of $A$
listed according to general multiplicity (see the description immediately before
Lemma~\ref{lem:ms}).
Notice that this is independent of the choice of $\lambda$ since
$E_\epsilon(\lambda \circ \sigma) = E_\epsilon(\lambda)$ for any permutation
$\sigma$ of $\{1,\ldots,k\}$.

\begin{theorem} 
Let $Z$ be a $\DT(\mu,c)$--operator, for any compactly supported Borel probability 
measure $\mu$ on the complex plane and any $c>0$.
Then $\delta_0(Z) =2$.
\end{theorem}

\begin{proof} Obviously $\delta_0(Z) \leq 2$ so it suffices to show
the reverse inequality.

We may without loss of generality assume $c=1$
(see~Proposition~2.12 of~\cite{dykema-haagerup:DT}).
Fix  $N \in \mathbb N$ with $N\ge2$.
By Theorem~4.12 of~\cite{dykema-haagerup:DT},
\begin{equation}\label{eq:Bmat}
\begin{bmatrix}
B_{11} & B_{12} & \cdots & B_{1N} \\ 0 & B_{22} & \ddots & \vdots \\ \vdots &
\ddots & \ddots & B_{N-1,N} \\ 0 & \cdots & 0 & B_{NN} \\ \end{bmatrix} \in\Mcal 
\otimes M_N(\mathbb C)
\end{equation}
is a $\DT(\mu,1)$--operator where $\{B_{11},\ldots, B_{NN} \} \cup
\langle B_{ij} \rangle_{1 \leq i < j \leq N}$ is a $*$-free family in $\Mcal$, the
$B_{ii}$ are $\DT(\mu, \frac{1}{\sqrt{N}})$--operators, and each $B_{ij}$ is
circular with $\varphi(|B_{ij}^2|) = \frac{1}{N}.$ From this we see that finding
microstates for $Z$ is equivalent to finding microstates for the operator~\eqref{eq:Bmat}
in $\Mcal\otimes M_N(\mathbb C)$.

Consider the sequence $\langle y_k \rangle_{k=1}^{\infty}$ constructed in Lemma 
3.2 and for each 
$\epsilon >0$ small enough, the corresponding sequence $\langle z_{k, \epsilon} 
\rangle_{k=1}^{\infty}$.
Let $R>1$, $m\in\Nats$, $\gamma>0$ and take $\gamma' = \gamma/16^m(R+1)^m>0$.
By Corollary~2.11 of \cite{dvv:strenghthened} there exist $k\times k$ complex unitary matrices $u_{1k},u_{2k},\ldots,u_{kk}$
such that $\{u_{1k}y_k u_{1k}^*,\ldots, u_{Nk} y_k u_{Nk}^* \}$ is an
$(m,\gamma')$--$*$--free family in $M_k(\mathbb C)$.
Also,by an application of Corollary~2.14 of~\cite{dvv:strenghthened}, there exists
a set $\Omega_k \subset \Gamma_R(\langle B_{ij} \rangle_{1
\leq i < j \leq N};m,k,\gamma')$ such that for any $\langle \eta_{ij} \rangle_{1
\leq i < j \leq N} \in \Omega_k$,
\[ \{u_{1k} y_k u_{1k}^*,\ldots, u_{Nk} y_k u_{Nk}^*\} 
\cup \langle \eta_{ij} \rangle_{1 \leq i < j 
\leq N}  
\]
is an $(m,\gamma')$-$*$ free family and such that
\begin{equation*}
\liminf_{k \rightarrow \infty} \left (k^{-2} \cdot \log
(\vol(\Omega_k)) + \frac{N(N-1)}{2} \cdot \log k \right) 
\ge 
\chi(\langle \RealPart B_{ij} \rangle_{1 \leq i < j \leq N},
\langle \ImagPart B_{ij} \rangle_{1 \leq i < j \leq N}) > - \infty,
\end{equation*}
where the volume is computed with respect to the product of the Euclidean norm $k^{1/2}|\cdot|_2$.
Since the operator~\eqref{eq:Bmat} is a
copy of $Z$, 
for any $\langle \eta_{ij}\rangle_{1 \leq i < j \leq N} \in \Omega_k$ we have
\[ \begin{bmatrix}
u_{1k} y_k u_{1k}^* & \eta_{12} & \cdots & \eta_{1N} \\ 0 & 
u_{2k}y_2 u_{2k}^* & \ddots 
&  \vdots \\ \vdots &
\ddots & \ddots & \eta_{N-1,N} \\ 0 & \cdots & 0 & u_{Nk} y_k 
u_{Nk}^*  \\ \end{bmatrix} \in \Gamma(Z;m, Nk,\gamma). \]

Because every complex matrix can be put into an upper-triangular form with
respect to an orthonormal basis, we can find for each $1 \leq j \leq N,$ a $k\times k$ unitary matrix
$v_{jk}$ such that $v_{jk} u_{jk} z_{k, \epsilon} u_{jk}^* v_{jk}^*$ is 
upper triangular. Observe now that for any
$\langle \eta_{ij} \rangle_{1\leq i < j \leq n} \in \Omega_k,$
\[  \begin{bmatrix}
v_{1k} & 0 & \cdots & 0 \\ 0 &
v_{2k} & \ddots
&  \vdots \\ \vdots &
\ddots & \ddots & 0 \\ 0 & \cdots & 0 & v_{Nk} \\ 
\end{bmatrix} 
\begin{bmatrix}
u_{1k} y_k u_{1k}^* & \eta_{12} & \cdots & \eta_{1N} \\ 0 &
u_{2k} y_k u_{2k}^* & \ddots
&  \vdots \\ \vdots &
\ddots & \ddots & \eta_{N-1,N} \\ 0 & \cdots & 0 & u_{Nk} y_k
u_{Nk}^*  \\ \end{bmatrix} 
\begin{bmatrix}
v_{1k}^* & 0  & \cdots & 0 \\ 0 &
v_{2k}^*  & \ddots
&  \vdots \\ \vdots &
\ddots & \ddots & 0 \\ 0 & \cdots & 0 & v_{Nk}^*.\\ \end{bmatrix} 
\]
is also an element of $\Gamma(Z;m,Nk,\gamma)$ and is 
equal to 
\[ \begin{bmatrix}
v_{1k} u_{1k} y_k u_{1k}^* v_{1k}^* & v_{1k} \eta_{12} 
v_{2k}^* & \cdots & v_{1j} \eta_{1N} v_{2k}^*\\ 0 &
v_{2j} u_{2k} y_k u_{2k}^* v_{2k}^* & \ddots
&  \vdots \\ \vdots &
\ddots & \ddots & v_{(N-1),k} \eta_{N-1,N} v_{Nk}^* \\ 0 & \cdots & 
0 & v_{Nk} u_{Nk} y_k
u_{Nk}^* v_{Nk}^* \\ \end{bmatrix}. \]

\noindent 
Moreover,
\[
|v_{jk}u_{jk}z_{k, \epsilon} u_{jk}^*v_{jk}^* - v_{jk}u_{jk} y_k 
u_{jk}^* v_{jk}^*|_2 = |z_{k, \epsilon} - y_k|_2
\]
and
$\limsup_{k\to\infty}|z_{k, \epsilon} - y_k|_2 \le \epsilon/\sqrt N$.
Therefore, for $k$ sufficiently large and for each $1 \leq j \leq N$  we have
$|v_{jk}u_{jk}z_{k, \epsilon} u_{jk}^*v_{jk}^* - v_{jk}u_{jk} y_k 
u_{jk}^* v_{jk}^*|_2\le\eps$.
Set $d_{jk} = v_{jk}u_{jk} z_{k,\epsilon} u_{jk}^* 
v_{jk}^*,$ and denote by $G_k$ the set of all $Nk \times Nk$ 
matrices of the form 

\[ \begin{bmatrix}
d_{1k} & v_{1k} \eta_{12}
v_{2k}^* & \cdots & v_{1j} \eta_{1N} v_{Nk}^*\\ 0 &
d_{2k} & \ddots
&  \vdots \\ \vdots &
\ddots & \ddots & v_{(N-1),k} \eta_{N-1,N} v_{Nk}^* \\ 0 & \cdots &
0 & d_{Nk} \\ \end{bmatrix}
\]

\noindent where $\langle \eta_{ij} \rangle_{1\leq i < j\leq N} \in
\Omega_k.$ Notice that each $d_{jk}$ is upper triangular and its 
eigenvalue distribution 
is exactly the same as that of $z_{k, \epsilon}$.
For $k$ sufficiently large, the set $G_k$ lies in the
$\epsilon$-neighborhood of $\Gamma(Z;m,Nk,\gamma)$.
Let $\theta(G_k)$ denote the unitary orbit of $G_k$ in $M_{Nk}(\Cpx)$.
We will now
find lower bounds for the $\epsilon$-packing numbers of 
$\theta(G_k)$ and 
thus, ones for $\Gamma(Z;m,Nk,\gamma).$

  Denote by $H_k \subset M_{Nk}(\mathbb C)$ all matrices of the 
form

\[ \begin{bmatrix}
0 & v_{1k} \eta_{12}
v_{2k}^* & \cdots & v_{1j} \eta_{1N} v_{Nk}^*\\ 0 &
0 & \ddots
&  \vdots \\ \vdots &
\ddots & \ddots & v_{(N-1),k} \eta_{N-1,N} v_{Nk}^* \\ 0 & \cdots &
\cdots &0 \\ \end{bmatrix}\]

\noindent where $\langle \eta_{ij} \rangle_{1\leq i < j \leq N} \in 
\Omega_k.$  Notice that $H_k$ is isometric to the space of 
all matrices of the form 

\[ \begin{bmatrix}
0 & \eta_{12} & \cdots & \eta_{1N} \\ 0 &
0 & \ddots
&  \vdots \\ \vdots &
\ddots & \ddots & \eta_{N-1,N}  \\ 0 & 
\cdot & \cdots &0 \\ \end{bmatrix}\]

\noindent where $\langle \eta_{ij} \rangle_{1 \leq i < j
\leq N} \in \Omega_k.$ It follows that $H_k$ must also have
the same volume as the above subspace, computed in the
obvious ambient Hilbert space of block upper triangular
matrices obeying the above decomposition.
Recall that for $n\in\Nats$, $T_n(\Cpx)$ denotes the set of uppertriangular matrices in $M_n(\Cpx)$;
let $T_{n,<}(\Cpx)$ denote the matrices in $T_n(\Cpx)$
that have zero diagonal, i.e.\ the strictly upper triangular
matrices in $M_n(\Cpx)$.
Denote by $W_k$
the subset of $T_{Nk,<}(\Cpx)$ consisting of all matrices $x$ such
that $|x|_2 < \epsilon$ and $x_{ij}=0$ whenever $1 \leq p < q \leq N$ and
$(p-1)k < i \leq pk$ and $(q-1)k < j \leq qk$.
Thus, $W_k$ consists of $N\times N$ diagonal matrices whose diagonal entries
are strictly upper triangular $k\times k$ matrices.
Denote by $D_k$ the subset of diagonal matrices $x$ of
$M_{Nk}(\mathbb C)$ such that $|x|_2 < \epsilon\sqrt2$. It
follows that if $f_k$ is the matrix

\[ \begin{bmatrix}
d_{1k} & 0 & \cdots & 0 \\ 0 &
d_{2k} & \ddots
&  \vdots \\ \vdots &
\ddots & \ddots & 0 \\ 0 & \cdots &
\cdots & d_{Nk} \\ \end{bmatrix}\]

\noindent then $f_k + D_k + W_k + H_k \subset \mathcal 
N_{3\epsilon}(G_k)$, where the $3\epsilon$ neighborhood is taken in
the ambient space $T_{Nk}(\Cpx)$ with respect to the metric induced by $|\cdot|_2$.
Now observe that the space of diagonal 
$Nk \times Nk$ matrices and $T_{Nk, <}(\Cpx)$ are orthogonal subspaces of $T_{Nk}(\Cpx)$.
Let $\theta_{3\eps}(G_k)$ denote the $3\eps$--neighborhood of the unitary orbit $\theta(G_k)$ of $G_k$.
Thus, denoting by $dX$ 
Lebesgue measure on $T_{Nk}(\Cpx)$ where $X = \langle x_{ij} \rangle_{1 \leq i \leq j 
\leq k}$, using Dyson's formula we have

\begin{eqnarray}
\notag \vol(\theta_{3 \epsilon}(G_k)) & \geq & C_{Nk} \cdot
\int_{f_k + D_k + W_k + H_k} \prod_{1 \leq i<j \leq Nk} |x_{ii} - x_{jj}|^2 dX \\[1ex]
\notag &= & C_{Nk} \cdot \vol(W_k + H_k) \cdot \int_{f_k + D_k} \prod_{1 \leq
i<j \leq Nk} |x_{ii} - x_{jj}|^2 dx_{11} \cdots dx_{(Nk)(Nk)} \\[1ex]
& \geq & C_{Nk} 
\cdot \vol(W_k + H_k) \cdot E_{\epsilon}(z_{k, \epsilon} \otimes I_N), \label{eq:NGamvol3}
\end{eqnarray}
where the constant $C_{Nk}$ is as in~\ref{eq:Ck}
and where $\vol(\theta_{3 \epsilon}(G_k))$ is computed in $M_{Nk}(\Cpx)$
and $\vol(W_k +H_k)$ is computed in $T_{Nk, <}(\Cpx)$, both being Euclidean volumes
corresponding to the norms $(Nk)^{1/2}|\cdot|_2$.
Clearly $\theta_{3\epsilon}(G_k) \subset \mathcal N_{4 
\epsilon}(\Gamma(Z;m,Nk,\gamma))$, so~\eqref{eq:NGamvol3} gives
a lower bound on $\vol(\mathcal N_{4 \epsilon}(\Gamma(Z;m,Nk,\gamma))$
as well.

Using~\eqref{eq:NGamvol3} and the standard volume comparison test, we have
\begin{eqnarray*}
P_{\epsilon}(\Gamma(Z;m,Nk,\gamma))
&\ge&\frac{\vol(\Nc_{4\eps}(\Gamma(Z;m,Nk,\gamma)))}{\vol(\Bc_{6\eps})} \\
& \ge & C_{Nk}\cdot E_{\epsilon}(z_{k,\epsilon}\otimes I_N)
\cdot\vol(W_k+H_k)
\cdot\frac{\Gamma((Nk)^2+1)}{\pi^{(Nk)^2} (6(Nk)^{1/2}\epsilon)^{2(Nk)^2}},
\end{eqnarray*}
where $\Bc_{6\eps}$ is a ball in $M_{Nk}(\Cpx)$ of radius $6\eps$ with respect to 
$|\cdot|_2$,
and we are computing volumes corresponding to the Euclidean norm $(Nk)^{1/2}|\cdot|_2$.
Since $W_k$ and $H_k$ are orthogonal, we have $\vol(W_k+H_k)=\vol(W_k)\vol(H_k)$,
where each volume is taken in the subspace of appropriate dimension.
But $W_k$ is a ball of radius $(Nk)^{1/2}\eps$ in a space of real dimension $Nk(k-1)$,
so
\[
\vol(W_k+H_k)=\frac{\pi^{\frac{Nk(k-1)}2}((Nk)^{1/2}\eps)^{Nk(k-1)}}{\Gamma(\frac{Nk(k-1)}2+1)}
\cdot (N^{1/2})^{k^2N(N-1)}\vol(\Omega_k).
\]
Applying Stirling's formula, we find
\begin{eqnarray*}  
\mathbb P_{\epsilon}(Z;m,\gamma) & \geq & \liminf_{k\to\infty}(Nk)^{-2}\log P_\eps(\Gamma(Z;m,Nk,\gamma)) \\[1ex]
& \geq & \liminf_{k\to\infty} (Nk)^{-2}\log(E_{\epsilon}(z_{k,\epsilon}\otimes I_N)) \\[0.5ex]
&& + \liminf_{k\to\infty}\bigg( \begin{aligned}[t]&(Nk)^{-2}\log(C_{Nk})+\frac1{2N}\log k + \frac1N\log\eps
 - \frac1{2N}\log(\frac{Nk(k-1)}2) \\[0.5ex]
    &+\log((Nk)^2)-\log k-2\log\eps+(Nk)^{-2}\log(\vol(\Omega_k))\bigg)+K_1 \end{aligned} \\[1ex]
&=& \liminf_{k\to\infty} (Nk)^{-2}\log(E_{\epsilon}(z_{k,\epsilon}\otimes I_N))
 + \liminf_{k\to\infty}\bigg((Nk)^{-2}\log(C_{Nk})+\frac12\log k\bigg) \\[0.5ex]
&& + \liminf_{k\to\infty}\bigg((Nk)^{-2}\log(\vol(\Omega_k))+(\frac12-\frac1{2N})\log k\bigg)
 + (2-N^{-1})|\log\eps|+K_2 \\[1ex]
&=& \liminf_{k\to\infty} (Nk)^{-2}\log(E_{\epsilon}(z_{k,\epsilon}\otimes I_N))
+N^{-2}\chi(\langle \RealPart B_{ij} \rangle_{1 \leq i < j \leq N},
\langle \ImagPart B_{ij} \rangle_{1 \leq i < j \leq N}) \\[0.5ex]
&& +\,(2-N^{-1})|\log\eps|+K_3,
\end{eqnarray*}
where $K_1$, $K_2$ and $K_3$ are constants independent of $\eps$, $m$ and $\gamma$.
Taking $m\to\infty$ and $\gamma\to0$, we get
\begin{eqnarray*}
\mathbb P_{\epsilon}(Z)&\ge&\liminf_{k\to\infty} (Nk)^{-2}\log(E_{\epsilon}(z_{k,\epsilon}\otimes I_N))
+N^{-2}\chi(\langle \RealPart B_{ij} \rangle_{1 \leq i < j \leq N},
\langle \ImagPart B_{ij} \rangle_{1 \leq i < j \leq N}) \\[0.5ex]
&&+\,(2-N^{-1})|\log\eps|+K_3.
\end{eqnarray*}
Since the eigenvalue distribution of $z_{k,\eps}\otimes I_N$ converges as $k\to\infty$ to the measure $\sigma_\eps$
of Lemma~\ref{lem:Brown},
dividing by $|\log \epsilon|$ and applying Lemma~\ref{lem:Eeps} now yields
\[
\delta_0(Z)=\limsup_{\epsilon \rightarrow0} 
\frac{\mathbb P_{\epsilon}(Z)}{|\log \epsilon|}
\ge\liminf_{\eps\to0}\frac{f(\epsilon)}{|\log \epsilon|} + 2- N^{-1} \ge 2-N^{-1}. 
\]
Since $N$ was arbitrary, it follows that $\delta_0(Z) \geq 2$, thereby 
completing the proof.
\end{proof}

\noindent{\it Acknowledgement.} A significant part of this research was 
conducted during the 2004 Free Probabilty Workshop at the Banff International Research Station,
and the authors would 
like to thank the organizers and sponsors for providing them with that opportunity to 
work together.
K.D.\ would like to thank the Mathematics Institute at the
Westf\"alische Wilhelms--Universit\"at M\"unster
for its kind hospitality during much of the time he was working on this project.

\end{document}